\documentclass{article}

\usepackage{graphicx}
\RequirePackage{url}
\usepackage{fancyvrb} 
\usepackage{amssymb,amsmath,amsfonts}
\usepackage{enumitem}
\usepackage{bm}
\usepackage{optidef}
\usepackage{empheq}

\newcommand{\ds}{\displaystyle}

\usepackage{soul}

\newcommand{\R}[1]{\mathbb{R}^{#1}}
\usepackage[citecolor=red,colorlinks=true,linkcolor=blue,linktoc=page]{hyperref} 

\usepackage{arxiv}

\usepackage[utf8]{inputenc} 
\usepackage[T1]{fontenc}    
\usepackage{hyperref}       
\usepackage{url}            
\usepackage{booktabs}       
\usepackage{amsfonts}       
\usepackage{nicefrac}       
\usepackage{microtype}      
\usepackage{lipsum}
\usepackage{graphicx}
\graphicspath{ {./images/} }

\title{Minimum Effort Control Using Variational Methods of Analytical Mechanics\\A New Approach For Optimal Control}

\author{
 Ossama Abdelkhalik \\
  Department of Aerospace Engineering\\
  Iowa State University\\
  Ames, Iowa 50011 \\
  \texttt{ossama@iastate.edu} \\
   \And
 Aimar Negrete \\
  Department of Aerospace Engineering\\
  Iowa State University\\
  Ames, Iowa 50011 \\
  \texttt{aimar@iastate.edu} \\
}

\begin{document}
\maketitle
\begin{abstract}
Modern optimal control theory involves adjoining the {already known} equations of motion of a dynamic system to the objective function using {dynamic costates}; this is done in order to constrain the optimal control solutions to satisfy the equations of motion. The use of costates increases the number of variables and hence increases the complexity of the problem. On the other hand, variational methods of analytical mechanics finds the equations of motion by minimizing an action functional of the dynamic system, realizing control forces as external input to the system. 
In this paper a new disruptive approach for computing the optimal control is presented. This approach adopts the variational methods of  analytical mechanics to derive equations for the control, in addition to the equations of motion. This is achieved by recognizing the control actuator as part of the dynamic system. In addition to the kinetic energy and potential energy, the action functional in this new approach includes additional energy terms that represent the control energy of the system. Two different methods are presented to write the modified action functional. The proposed approach is a significant departure from the modern optimal control theory, and it eliminates the need for costates when solving for the control. In this paper, a case study is presented to demonstrate the new approach.
\end{abstract}


\section{Introduction}

The ability to solve Optimal Control Problems (OCPs)  in a computationally efficient way is necessary for numerous engineering applications. There are several methods developed over the past several decades to solve OCPs that can be broadly categorized into Direct Methods, Indirect Methods, and Dynamic Programming\cite{Bohme2017} . 
In a direct method, the state and/or control are approximated and the OCP is transcribed to a nonlinear optimization problem (or nonlinear programming problem) \cite{AnilRao2009Survey}. 
In a wide range of applications, dynamic programming breaks complex problems into simpler sub-problems that are solved recursively as detailed in many references \cite{Bertsekas1987DP}. The indirect approach provides necessary conditions for optimality, but finding a solution that satisfies these conditions for problems with a significant number of controls, nonlinear models, or distributed parameter systems can be highly challenging. Despite these difficulties, indirect methods will outperform any other solution method in terms of accuracy \cite{BohmeThomas2017}.

The classical optimal control theory is a concomitant of the calculus of variations, a branch of mathematical analysis concerned with finding the stationary point of a functional \cite{variations}. Roots for the OCP can also be found in classical control theory, and in linear and nonlinear programming \cite{Bryson1996Article}. 
There is indeed a wealth of literature that details the modern optimal control theory and how to apply it \cite{Bryson1975}. 
The basic process of solving an OCP using the indirect approach involves adjoining the equations of motion of a dynamic system to the objective function using \emph{dynamic costates}; this is done in order to constrain the OCP solutions to satisfy the system's equations of motion. Specifically, a costate variable would be introduced for each state variable in the corresponding system. Consider for example an OCP of the form:
\begin{eqnarray} \label{OCPformulation}
   \textbf{Minimize} & \, \, J =  & \Phi({\bm{q}}(t_f),t_f)+\int_{t_0}^{t_f}L({\bm{q}},{\bm{u}},t)dt \\  
    \textbf{Subject to} &  \, \, \dot{{\bm{q}}} = & {f}({\bm{x}},{\bm{u}},t) \label{OCPformulation1}
\end{eqnarray}
where ${\bm{q}}$ is the state vector, ${\bm{u}}$ is the control vector, $t$ is time, $\Phi$ represents cost from and constraints on the state at the terminal time $t_f$, and the continuous (or incremental) cost is represented by $L$. The constraint in \eqref{OCPformulation1} is the equations of motion.
Classical optimal control theory introduces the Hamiltonian function for this system as $H = L + \bm{\lambda}^T\bm{f}$, where $\bm{\lambda}$ are {the dynamic costates appended to the system} to enforce the differential equation constraints. The necessary conditions for optimality are 
\begin{equation}
\begin{split}
    \dot{\bm{\lambda}}^T = -\frac{\partial H}{\partial \bm{q}} \, \, , \,  \,
    \frac{\partial H}{\partial\bm{u}} = 0  \, \, , \,  \,
    \bm{\lambda}(t_f) = \frac{\partial \Phi}{\partial \bm{q}(t_f)}
\end{split}
\label{ELElambda}
\end{equation}

Through the Euler-Lagrange equations \cite{variations}, differential equations and terminal conditions are derived for the costates in terms of, in general, the states, control, and time. The result is a two-point boundary value problem in the states and costates that - when solved - returns the optimal trajectories of the system. The Legendre-Clebsch condition determines whether the stationary solution is a minimum or a maximum. Pontryagin's Minimum Principle applies for the more general problem with inequality constraints. 

There are several challenges when implementing this current approach for optimal control; one challenge for instance, in many classes of problems, the terminal states are constrained. This shows up for example in space trajectory optimization problems, where it is desired to match the target position and velocity vectors at the final time, whether it be for a planetary rendezvous or orbit insertion. When this occurs, no information is known at the endpoints about the costates, as each of the terminal conditions on the costates is now a Lagrange multiplier to be solved for. Finding a good initial guess for the costates, and subsequently a good solution, is not a trivial task.

The new approach for optimal control presented in this paper, however, is a significant departure from the modern optimal control theory referred to above. It diverges from existing methods in that it does not require knowledge of the equations of motion; the equations of motion are outcomes along with the equations for optimal control when using this new approach. Moreover, the new approach does not use costates in computing the optimal control. This new approach is rooted in the elegant Principle of Least Action as described in Section~\ref{Hypo}. The next paragraph briefs necessary background on variational methods of analytics mechanics.

At the same time as the calculus of variations was being developed, several new advancements in analytical mechanics were established that culminated into the principle of least action. In general, these advancements are concerned with finding paths that satisfy $\delta S = 0$, where $\delta(\cdot)$ is the variational operator \cite{variations} and $S$ is a function that describes the action of the system. For clarity, the Lagrangian $\mathcal{L} \equiv T-V$ is defined as the difference in kinetic and potential energies of a system in terms of its generalized coordinates $\bm{q}$ \cite{AM-Schaub-2nd}. The original principle of least action comes from Pierre Louis Maupertuis, which gives the action integral as $S_M=\int_a^b \bm{p}\cdot\bm{dq}$, where $p_i \equiv \partial L/\partial \dot{q}_i$ are the conjugate momenta of each generalized coordinate. Finding the stationary point of this action integral yields the trajectory with respect to the generalized coordinates. This approach returns the path of the trajectory, rather than equations of motion. A similar approach is Hamilton’s law of varying action \cite{mechanics}, which instead defines the action integral as $S_H = \int_a^b (T+W)dt$, where $W$ is the work done on the system. The stationary point of this action yields the trajectories with respect to time. This returns the equations of motion of a system directly. There exist  other formulations along the same lines, such as Lagrangian dynamics (a special case of Hamilton's law), the Gibbs-Appell equations of motion, and Gauss’ principle of minimum constraint, among others. These formulations vary in both the underlying mathematics and the types of equations returned by the process. 
With the variational methods of analytical mechanics, the equations of motion of a system are derived by finding the stationary point of an action functional. Under necessary second order conditions, this is equivalent to minimizing the action functional. Optimal control theory by parallel aims to find the optimal input to a dynamic system that minimizes some objective functional. It is possible to write the relationship between the maximum principle, which is the fundamental idea of optimal control theory, and classical mechanics as detailed in \cite{Todorov2006OCP}.

This paper presents the main concept of the proposed approach, and presents several detailed examples highlighting the process for computing the optimal control. Section~\ref{Hypo} of this paper presents the hypothesis of the new approach.
Section~\ref{ProposApp} presents the proposed approach. We present two methods for finding the variational principle, and in each case provide a case study. Section ~\ref{invariancesSection} presents how to find new invariances for the optimal control system using the variational principles found in Section ~\ref{ProposApp}, a useful technique for solving optimal control problems once the variational principle is found. 

\section{The Hypothesis}\label{Hypo}
In the classical approach for solving the OCPs, the control generalized force is considered as external input to the system; the equations of motion of such system can be obtained using, for instance, Hamilton's Principle. According to Hamilton's Principle, the equations of motion are obtained by minimizing an action integral that is usually written in terms of the kinetic energy and potential energy of the system. This means that the system states will follow a trajectory that minimizes the action integral. 

In the new approach, we consider the actuator that provides the control force as an integral part of the system; hence the control force is now a state. When writing the equations of the motion for the new system using Hamilton's Principle, the equations should include also equations for the control (since it is now a state). We also modify the action integral to account for the energy available for the actuator (here referred to as the control energy); for example if we are using a chemical engine onboard a spacecraft then the energy available for the actuator would be the chemical energy content in the amount of fuel available onboard. Hence, the new action integral is a function of the kinetic energy, the potential energy, and the control energy. 

The hypothesis of this approach is as follows: Since equations of motion follow by minimizing an action integral that is a function of kinetic energy and potential energy, then minimizing the new action integral that is a function also of the control energy would yield equations of motion as well as equations for the minimum effort control. 

The rest of this paper presents a complete derivation and demonstration by example for this new approach.

\section{The New Optimal Control Approach}\label{ProposApp}
\begin{figure}[hbt!]
    \centering
    \includegraphics[width=0.7\textwidth]{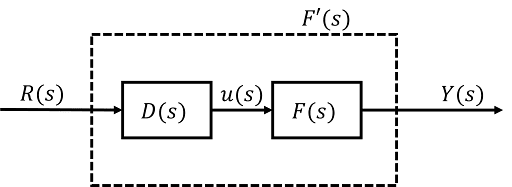}
    \caption{Block diagram for a typical dynamic system}
    \label{BlockDiag}
\end{figure}

In modern control theory, a dynamic system is controlled via an actuating system. In Fig.~\ref{BlockDiag}, the system dynamics are represented by $F(s)$ and the actuator dynamics are represented by $D(s)$. A reference signal $R(s)$ is passed into the overall system, and once the signal has passed through both the actuator dynamics and the system dynamics, the output from the complete system is $Y(s)$.

Classical optimal control theory is concerned with minimizing a functional of the form given in Eq.~\eqref{OCPformulation} subjected to differential equation constraints. 
When applying classical optimal control theory to a problem, it is necessary that the governing differential equations are known a priori. These can be derived using Hamilton's law of varying action. 
{To write the equations of motion of a holonomic system with Hamilton's law of varying action, one begins by defining the action functional}
\begin{equation*}
S_H=\int_{t_0}^{t_f}(T+W)\,dt
\end{equation*}
where $T$ represents the system's total kinetic energy, and $W$ is the sum of conservative and nonconservative work done on the system. The integrand of the action functional is known as Hamilton's principal function \cite{AM-Schaub-2nd}. Expressing the total work as $W = W_{nc} - V$, the action functional can be rewritten as 
\begin{equation*}
S_H=\int_{t_0}^{t_f}(\mathcal L + W_{nc})\,dt
\end{equation*}
where $\mathcal L = T - V$ is the traditional Lagrangian of the dynamic system as defined previously. Next, the variational methods of analytical mechanics can be used to derive the equations of motion (or equivalently $F(s)$ in Figure~\ref{BlockDiag}) by {rendering stationary} the action integral $S_H$, {under the assumption that the variations of the system states vanish at the endpoints.} These equations of motion include the control signal $u(t)$ as an external input to the system, contained in the $W_{nc}$ term. In the proposed approach, the control actuation is considered part of the dynamic system, and the same variational methods of analytical mechanics are used to write equations for this new system (represented as $F'(s)$ in Figure~\ref{BlockDiag}.) The resulting equations will be shown here to include the equations of motion of the original system (e.g. $F(s)$) as well as equations for the control $u(t)$. This is achieved by  a modified action integral $S'=\int L\,dt$. Note that the ``optimal control Lagrangian", $L$, is distinct from the classical Lagrangian, $\mathcal{L}$. In this paper, we are concerned with methods for finding the Lagrangian $L$.

The main difference in the proposed approach when compared to the modern optimal control theory is that the control, rather than being treated as an input to the system, is treated as an additional state in the system. As such, the action functional $S_H$ must be rewritten in such a way as to include also the ``control energy". By then minimizing the modified action functional as in the traditional methods of analytical mechanics, treating the control as a state, the original equations of motion of the system are returned, in addition to new equations of motion which describe the optimal evolution of the control. To do this, a method is developed in section~\ref{InvarHamilt} that leverages the fact that the optimal control Hamiltonian is invariant on the optimal trajectory to derive the Lagrangian; this method involve a total transformation of the action functional into a form which couples directly the action of the states and the control through the application of Noether's theorem {Specifically, section \ref{InvarHamilt} treats the case of conservative holonomic physical systems.}. In section~\ref{SecInvLagrang} a method is developed to write the Lagrangian within the framework of the Helmholtz conditions; this approach assumes knowledge of the optimal control, and develops expressions of the Lagrangian. The Helmholtz condition method is more general than the Hamiltonian invariance method. {This more general case allows for the treatment of nonconservative systems. In both cases, the systems considered are finite dimensional.} Regardless of which method is used to find $L$, the Lagrangian is then used to write motion invariants on the optimal control trajectory in section~\ref{invariancesSection}.

\subsection{Invariance of The Hamiltonian}\label{InvarHamilt}
For any conservative system, such as the simple mass spring system with no control force, the total energy $E=T+V$ is a conserved quantity of the system. The difference between the kinetic and potential energies is the classical Lagrangian, $\mathcal{L}$, whose time integration is the action integral to be minimized to get the equations of motion in the traditional analytical mechanics approach. In the classical optimal control theory, the Hamiltonian function is a conserved quantity of the system on the optimal path when the equations of motion and the continuous cost are not explicit functions of time. The optimal control problem now becomes to find a Lagrangian function whose system equations of motion render the Hamiltonian invariant.

To do this, the mathematical framework considered is Noether's theorem, a connection between symmetries of a system and corresponding conserved quantities. 
Consider the arbitrary functional
\begin{equation}
\ds \Gamma=\int_a^b L(t,\bm{q},\dot{\bm{q}}) dt
    \label{NoetherThermActInteg}
\end{equation}
Consider also the infinitesimal time and generalized coordinates transformations, respectively:
\begin{equation}
\ds t'=t+\varepsilon\xi(t,\bm{q})+\mathcal{O}(\varepsilon^2)
    \label{NoetherThermTimeTrans}
\end{equation}
\begin{equation}
\ds \bm{q}'=\bm{q}+\varepsilon\bm{\eta}(t,\bm{q})+\mathcal{O}(\varepsilon^2)
    \label{NoetherThermCoordTrans}
\end{equation}
where $\varepsilon << 1$. 

The integral in \eqref{NoetherThermActInteg} is said to be invariant if under the time and generalized coordinate transformations \eqref{NoetherThermTimeTrans} and \eqref{NoetherThermCoordTrans} the value of the integral remains arbitrarily close to the original. That is, the integral is invariant if and only if 
\begin{equation*}
\ds \left|\int_a^b L(t,\bm{q},\dot{\bm{q}}) dt - \int_a^b L(t',\bm{q}',\dot{\bm{q}}') dt\right| \sim  \varepsilon^s
    \label{NoetherTherm01}
\end{equation*}
for $s>1$. Since these transformations are infinitesimal, the integrals are only required to be invariant to the linear term, as in the limit as $\varepsilon\rightarrow 0$, they will be identical.
Noether's theorem shows that for each transformation $\xi,\bm{\eta}$ (in \eqref{NoetherThermTimeTrans} and \eqref{NoetherThermCoordTrans}) that renders the action integral stationary, there will be a corresponding conserved quantity of the system given by
\begin{equation*}
\frac{\partial L}{\partial \dot{\bm{q}}}(\bm{\eta}-\xi\dot{\bm{q}}) + L\xi\equiv\mathrm{const}.
    \label{NoetherThermInv}
\end{equation*}
For a conservative physical system, the Lagrangian will be invariant under a transformation of $\xi=1$, $\bm{\eta} = 0$, known as a time transformation. This leads to a conservation law which corresponds to a conservation of the total energy of the system, $T + V = \mathrm{const.}$ The idea behind this approach here is that the Hamiltonian is a conserved quantity in the unconstrained optimal control problem for autonomous systems. If we treat the optimal control Hamiltonian as a total ``energy" of the optimal control system, then we can say that $H = T'+V'$, where $T'$ and $V'$ are ``kinetic" and ``potential" energies of the optimal control system in terms of state and control variables. Under a time transformation, this corresponds to an optimal control Lagrangian $L = T' - V'$. One can therefore arrive at a Lagrangian for an optimal control system by first rewriting the Hamiltonian as a function of only the state and control by eliminating the costates, then flipping the sign on the strictly position dependent terms. This approach follows directly from the idea that the Hamiltonian is the conservation law for a time transformation, allowing one to work backwards to arrive at the Lagrangian. As can be observed this approach requires writing the Hamiltonian, hence the dynamic costates. However, in the following section, we provide a method for arriving at the optimal control Lagrangian with knowledge only of the original physical system.  
Without loss of generality, the function $f(q,\dot{q})$ in Eq. \eqref{OCPformulation1} can be written in the form
\begin{equation*}
\bm f=\bm g(\bm q,\dot{\bm q}) + \bm h(\bm q)
\end{equation*}
If we conduct the process described above, we arrive at the following general form for the Lagrangian
\begin{equation}
{L} = \dot{\bm u}^T\dot{\bm q} + \bm u^T\frac{\partial \bm f}{\partial \dot{\bm{q}}} \dot{\bm q} - \bm u^T\bm g + \bm u^T \bm h + \frac{1}{2}\bm u^T\bm u
    \label{BetaFunctional}
\end{equation}
Writing this Lagrangian for an arbitrary physical system requires no knowledge of the optimal control solution or any costate information - not even the form of the Hamiltonian. Further along, it will even be possible to eliminate $f, g,$ and $h$ from $L$, writing them instead in terms of the original system's classical Lagrangian $\mathcal{L}$. 
The following theorems prove that the Lagrangian given in \eqref{BetaFunctional} yields an equation for the optimal control when minimized with respect to the generalized coordinates, and the original state differential equations when minimized with respect to the control variables.

First we use the classical optimal control theory to prove a differential equation form for the following optimal control problem. Suppose that the system equations of motion using the generalized coordinates \cite{AM-Schaub-2nd} take the form: 
\begin{equation}
  \ddot{\bm q} = \bm f(\bm q,\dot{\bm q},t) + \bm u
    \label{SysDyn}
\end{equation}

\noindent The minimum energy optimal control problem can be formulated as

\noindent\textbf{Problem 1:} \label{OCProbelm}\begin{mini}|s| 
{\mathbf{x}} { \mathcal{J} = \bm \nu^T\bm{F}(\bm{x}(t_f),t_f) + \frac{1}{2}\int_{t_0}^{t_f}\bm u^T\bm u\,dt}{}{}
\addConstraint{ \ddot{\bm q} = \bm f(\bm q,\dot{\bm q},t) + \bm u}
\end{mini}where the terminal constraints vector is $\bm F\in\R{p}$ and $\bm \nu\in\R{p}$ is the associated vector of Lagrange multipliers.

\textbf{Theorem}
    For the system defined in \eqref{SysDyn}, the optimal control solution to the minimum control effort problem defined in \textbf{Problem 1} will satisfy $$\ddot{\bm u} = -\left(\frac{\partial \bm f}{\partial \dot{\bm q}}\right)^T\dot{\bm u} - \left[ \frac{d}{dt}\left(\frac{\partial \bm f}{\partial \dot{\bm q}} \right) - \frac{\partial \bm f}{\partial \bm q} \right]^T\bm u$$

\textbf{proof}
Let the state vector $\bm{x}\in\R{2n}$ be defined as $\bm x=\left[ \begin{array}{c} \bm q \\ \hline \dot{\bm q} \end{array} \right]$. Additionally, let the costate vector $\bm \lambda\in\R{2n}$ be defined as $\bm \lambda=\left[ \begin{array}{c}\bm \Lambda_1 \\ \hline \bm \Lambda_2 \end{array} \right]$, where $\bm \Lambda_1,\bm \Lambda_2\in\R{n}$ are the costate vectors associated with $\bm q$ and $\dot{\bm q}$, respectively. With these definitions, the optimal control Hamiltonian function for the problem is given by 

\begin{equation}
        H = \frac{1}{2}\sum u_i^2 + \bm \Lambda_1^T \bm q + \bm \Lambda_2^T(\bm f + \bm u)
    \label{hamiltonianOCP}
\end{equation}

\noindent From the necessary conditions of optimality \cite{Bryson1975}, the costates will have differential equations given by $\dot{\lambda}_i = -\frac{\partial H}{\partial x_i}$. For our system,

\begin{align}
    \dot{\bm\Lambda}_1 &= -\left(\frac{\partial \bm f}{\partial \bm q}\right)^T\bm \Lambda_2 \label{lambda1dot}\\ 
    \dot{\bm \Lambda}_2 &= -\bm \Lambda_1 - \left(\frac{\partial \bm f}{\partial \dot{\bm q}}\right)^T\bm \Lambda_2 \label{lambda2dot}
\end{align}

\noindent The stationarity condition requires that $\frac{\partial H}{\partial \bm u} = 0$. For this system, $\frac{\partial H}{\partial \bm u} = \bm u^T + \bm \Lambda_2^T = \bm 0$. This gives the relations 

\begin{equation}
    \begin{split}
        \bm u = -\bm \Lambda_2 \, , \,   \dot{\bm u} = -\dot{\bm \Lambda}_2 \, , \, 
        \ddot{\bm u} = -\ddot{\bm \Lambda}_2
    \end{split}
    \label{dhdu}
\end{equation}

\noindent Passing in the relations in Eq. \eqref{dhdu} into Eqs. \eqref{lambda1dot} and \eqref{lambda2dot} yields 

\begin{align}
    \dot{\bm \Lambda}_1 &= \left(\frac{\partial \bm f}{\partial \bm q} \right)^T \bm u \label{intermediate1} \\ 
    -\dot{\bm u} &= -\bm \Lambda_1 + \left( \frac{\partial \bm f}{\partial \dot{\bm q}} \right)^T \bm u \label{intermediate2}
\end{align}

\noindent By differentiating Eq. \eqref{intermediate2} with respect to time and substituting in Eq. \eqref{intermediate1}, we arrive at 

\begin{equation}
    \ddot{\bm u} = -\left(\frac{\partial \bm f}{\partial \dot{\bm q}} \right)^T \dot{\bm u} - \left(\frac{d}{dt}\left(\frac{\partial \bm f}{\partial \dot{\bm q}}\right) - \frac{\partial \bm f}{\partial \bm q}\right)^T\bm u
    \label{optimalControlODE}
\end{equation}
Equation~\eqref{optimalControlODE} is the set of differential equations that the minimum effort optimal control solution must satisfy.\setlength\fboxrule{1.1pt}\setlength\fboxsep{0pt}\fbox{\phantom{\rule{5pt}{5pt}}}



{
\textbf{Theorem}\label{TheoremLagrange}
    {(Euler-Lagrange Equation)} The functional $\int_a^b F(t,\bm y,\dot{\bm y}) \,dt$ is {rendered stationary} by the {family of functions} which {satisfy} $$\frac{d}{dt}\left(\frac{\partial F}{\partial \dot{\bm y}}\right) - \frac{\partial F}{\partial \bm y} = \bm 0$$
}

\textbf{Proof}
The variation of a functional $I = \int G(t,\bm x) dt$ is defined as \cite{gelfand2000calculus} $$\delta I = \int \sum \frac{\partial G}{\partial x_i}\delta x_i\,dt$$ where $\delta x_i(t)$ is the variation of coordinate $x_i$ at time $t$.  The functional will be minimized when $\delta I = 0$. Taking the variation of the functional in this theorem, 

\begin{equation}
    \delta \int_a^b F(t,\bm y,\dot{\bm y})\,dt = \int_a^b \left( \frac{\partial F}{\partial \dot{\bm y}}\delta \dot{\bm y} + \frac{\partial F }{\partial\bm  y}\delta \bm y \right) dt
    \label{deltaF1}
\end{equation}

\noindent The variations $\delta \dot{\bm y}$ and $\delta \bm y$ are related, so the equation must be integrated by parts. 

\begin{equation}
        \int_a^b \frac{\partial F}{\partial \dot{\bm y}} \delta \dot{\bm y}\,dt = \left.\frac{\partial F}{\partial \dot{\bm y} } \delta \bm y \right|_a^b - \int_a^b \frac{d}{dt}\left(\frac{\partial F}{\partial \dot{\bm y}} \right) \delta \bm y\,dt
        \label{deltaydot}
\end{equation}

\noindent For many applications, it is common to assume that the variations vanish at the endpoints; that is $\delta \bm y(a) = \delta \bm y(b) = \bm 0$ \cite{AM-Schaub-2nd}. With this assumption, we can combine Eqs. \eqref{deltaF1} and \eqref{deltaydot} to write the total variation of the functional.

\begin{equation*}
    \delta\int_a^b F(t,\bm y,\dot{\bm y})dt = -\int_a^b \left(\frac{d}{dt}\left( \frac{\partial F}{\partial \dot{\bm y}}\right) - \frac{\partial F}{\partial \bm y} \right) \delta \bm y \,dt
\end{equation*}

\noindent The minimum of this functional will be found when the variation vanishes. The quantity $\delta \bm y$ is arbitrary and in general nonzero, so the only way to ensure that the integral vanishes is for the bracketed term to be identically zero; that is to say 
\begin{equation*}
    \frac{d}{dt}\left(\frac{\partial F}{\partial \dot{\bm y}} \right) - \frac{\partial F}{\partial \bm y} = \bm 0
\end{equation*}\setlength\fboxrule{1.1pt}\setlength\fboxsep{0pt}\fbox{\phantom{\rule{5pt}{5pt}}}

To arrive at an optimal control functional, we apply Noether's theorem on the optimal control Hamiltonian function in Eq. \eqref{hamiltonianOCP} under the assumption that the Hamiltonian is in fact an invariant of the system under the time transformation $\xi=1,\bm{\eta}=0$. First, $H$ must be rewritten to exclude the costates from the function in favor of the control and its derivatives. Using the relations developed in Eqs. \eqref{dhdu} and \eqref{intermediate2}, we can rewrite the Hamiltonian as

\begin{equation}
    H = \frac12 \bm u^T \bm u + \left[ \dot{\bm u} + \left(\frac{\partial \bm f}{\partial \dot{\bm{q}}}\right)^T\bm u \right]^T \dot{ \bm{q}} - \bm{u}^T(\bm f + \bm u)
    \label{hamiltonian1}
\end{equation}

$\newline$
\noindent\textbf{Assumption 1} We assume that the system described in Equation~\eqref{SysDyn} is conservative and autonomous; hence the function $\bm f$ in Eq.~\eqref{SysDyn} can be expressed as $\bm f(\bm q,\dot{\bm q})=\bm g(\bm q,\dot{\bm q}) + \bm h(\bm q)$. 

With this assumption, Eq. \eqref{hamiltonian1} can be explicitly split into strictly position dependent terms, and velocity dependent terms. We note that the position terms in the extended generalized coordinates are functions of strictly $\bm{q}$ and $\bm{u}$, and the velocity dependent terms are all functions of $\dot{\bm{q}}$ and $\dot{\bm{u}}$.

\begin{equation*}
    H = \left[ \dot{\bm{u}}^T\dot{\bm q} + \bm{u}^T \frac{\partial \bm f}{\partial \dot{\bm q}}\dot{\bm q} - \bm{u}^T \bm g \right] + \left[-\bm{u}^T\bm h - \frac{1}{2}\bm{u}^T\bm{u} \right] = T' + V'
\end{equation*}

\noindent Treating the Hamiltonian now as a total energy for the optimal control system, and under the assumption that it is invariant under a time transformation ($\xi = 1, \bm\eta=\bm 0$), the Lagrangian for the optimal control system can be written simply as $L = T' - V'$. Next we use \textbf{Theorem \ref{TheoremLagrange}} to find the solution to the optimal control Lagrangian problem below.

\textbf{Theorem}\label{TheormOptSol}
For the dynamic system described by Eq.~\eqref{SysDyn} which satisfies \textbf{Assumption 1}, and considering the function
\begin{equation}\label{BetaFunctional2}
{L} = \dot{\bm u}^T\dot{\bm q} + \bm u^T\frac{\partial \bm f}{\partial \dot{\bm{q}}} \dot{\bm q} - \bm u^T\bm g + \bm u^T \bm h + \frac{1}{2}\bm u^T\bm u
\end{equation} where $\bm f = \bm g(\bm q,\dot{\bm q}) + \bm h(\bm q)$, finding the {stationary point} of the functional $\int_{t_0}^{t_f}L\,dt$ over the generalized coordinates $\bm q$

\begin{mini}|s| 
{\bm q} {\int_{t_0}^{t_f}L\,dt}{}{}\label{OCPAuxrobelm}
\end{mini}

\noindent will provide the optimal control differential equations, and {stationary} with respect to the control variables $\bm u$ will provide the original state differential equations. 

\textbf{Proof}
Let $\mathcal{L} = T-V$ be the classical Lagrangian function of the given system, where $T$ is the system's kinetic energy, and $V$ is the system's potential energy. The system equations of motion can be written in the form \cite{sarlet}:
\begin{equation}
    \frac{d}{dt}\left( \frac{\partial \mathcal{L}}{\partial \dot{\bm q}} \right) - \frac{\partial\mathcal{L}}{\partial \bm q} = \alpha(\ddot{\bm q} - \bm f)
    \label{lagrangianEquation}
\end{equation}
\noindent where $\ds \alpha = \ds \frac{\partial^2 \mathcal{L}}{\partial \dot{\bm q}^2}$ is known as the Jacobi Last Multiplier (JLM) matrix \cite{nucci2008jacobi}. The properties of this matrix will be further explored in Section \ref{SecInvLagrang}, but in {general} a proper JLM is a nonsingular symmetric matrix. Note that in the case of the uncontrolled conservative system, $\ddot{\bm q}=\bm f$, and hence the right hand side vanishes. 

For a conservative system, the position vector will be a function of only the generalized coordinates - that is $r_i=r_i(q)$. For a general physical system, the kinetic energy is of the form $T = \frac{1}{2}\sum m_i |\dot{r}_i|^2$. From the definition of $\bm r$,

\begin{align*}
    \dot{r}_i  = \frac{dr_i}{dt} = \sum_j\frac{\partial r_i}{\partial q_j} \frac{\partial q_j}{\partial t} = \sum_j \frac{\partial r_i}{\partial q_j} \dot{q}_j
\end{align*}

\noindent The kinetic energy is

\begin{equation}
\begin{split}
    T &= \frac{1}{2}\sum_i m_i \left| \sum_j \frac{\partial r_i}{q_j}\dot{q}_k \right|^2 \\ 
    &= \frac{1}{2}\sum_i m_i \sum_{j,k}\left( \frac{\partial r_i}{\partial q_k} \cdot \frac{\partial r_i}{\partial q_j} \right) \dot{q}_j \dot{q}_k \\
    &= \frac{1}{2}\sum_{j,k}M_{jk}\dot{q}_j\dot{q}_k
\end{split}
\label{kineticenergy}
\end{equation}

\noindent We can write the classical Lagrangian as 
\begin{equation}
    \mathcal{L} = T - V = \frac{1}{2}\dot{\bm q}^T M \dot{\bm q} - V(\bm q)
    \label{generalLagrangian}
\end{equation}

\noindent Note that in general, the mass matrix $M\in \R{n\times n}$ is positive semidefinite, symmetric, and invertible. From the definition of the Jacobi Last Multiplier matrix $\alpha$, it is clear that for this class of systems $\alpha = M$. We can use Eq. \eqref{lagrangianEquation} to find the equations of motion for the system from the Lagrangian function. Computing the partial derivatives of Eq. \eqref{generalLagrangian} with respect to $\bm q$ and $\dot{\bm q}$ yields

\begin{align}
    \frac{\partial \mathcal{L}}{\partial \bm q} &=  \frac{1}{2}\frac{\partial }{\partial \bm  q}\dot{\bm q}^T M \dot{\bm q} - \frac{\partial V}{\partial \bm q}\\ 
    \frac{\partial \mathcal{L}}{\partial \dot{\bm q} } &= \dot{\bm q}^T M
    \label{dldqold}
\end{align}

\noindent The total time derivative of $\partial \mathcal{L}  / \partial \dot{\bm q}$ is given by 

\begin{equation}
    \begin{split}
        \frac{d}{dt}\left(\frac{\partial \mathcal L}{\partial \dot{\bm q}}\right) &= \ddot{\bm q}^T M + \dot{\bm q}^T \dot M \\
        &= \ddot{\bm q}^T M + \dot{\bm q}^T\left( \sum_i \frac{\partial M}{\partial q_i}\frac{\partial q_i}{\partial t} + \frac{\partial M}{\partial t} \right) \\ 
        &= \ddot{\bm q}^T M + \dot{\bm q}^T \left( \sum_i \frac{\partial M}{\partial q_i}\dot{q}_i\right)
    \end{split}
    \label{ddtdldqold}
\end{equation}

\noindent From the definition of $M$ in Eq. \eqref{kineticenergy}, $M$ is only a function of the generalized coordinates $\bm q$, hence the above expansion of $\dot{M}$. 
Combining this with Eq. \eqref{ddtdldqold} into Eq. \eqref{lagrangianEquation} yields the final form for the equations of motion of the system.

\begin{equation}
    \bm f = M^{-1}\left[ \frac{1}{2}\dot{\bm q}^T\frac{\partial M}{\partial \bm  q} \dot{\bm q} - \frac{\partial V}{\partial \bm q} - \dot{\bm q}^T\left(\sum_i \frac{\partial M}{\partial q_i}\dot{q}_i\right) \right]
    \label{feq}
\end{equation}

\noindent From \textbf{Assumption 1}, the function $\bm f(\bm q,\dot{\bm q})$ can be written in the form
\begin{equation*}
\bm f=\bm g(\bm q,\dot{\bm q}) + \bm h(\bm q)
    \label{DynForm}
\end{equation*}
The $\bm g$ functions have the velocity dependent terms in the equations of motion, and $\bm h$ are the strictly position based terms. Therefore, $$\frac{\partial\bm  f}{\partial \dot{\bm q}} = \frac{\partial \bm g}{\partial \dot{\bm q}} + \frac{\partial \bm h }{\partial \dot{\bm q}} = \frac{\partial \bm g}{\partial \dot{\bm q}}$$ From Eq. \eqref{feq}, we can write that $$g=M^{-1}\left[ \frac{1}{2}\dot{\bm q}^T\frac{\partial M}{\partial \bm q} \dot{\bm q} - \dot{\bm q}^T\left(\sum_i \frac{\partial M}{\partial q_i}\dot{q}_i\right) \right]$$ Since this is quadratic in the generalized velocities, it will satisfy the relation 

\begin{equation}
    \frac{\partial \bm g}{\partial \dot{\bm q}}\dot{\bm q} = 2\bm g
    \label{quadraticSatisfaction}
\end{equation}


To {render stationary} the functional $\int L\,dt$ with respect to the generalized coordinates $\bm q$, we will use \textbf{Theorem \ref{TheoremLagrange}}. We first compute the partial derivatives of $ L$ with respect to $\bm q$ and $\dot {\bm q}$. 

\begin{align}
    \frac{\partial  L}{\partial \bm q} &= \bm u^T\left(\frac{\partial}{\partial \bm q}\frac{\partial \bm f}{\partial \dot{\bm q}}\right)\dot{\bm q} - \bm u^T \frac{\partial \bm g}{\partial \bm q} + \bm u^T \frac{\partial \bm h}{\partial \bm q} \\ 
    \frac{\partial  L}{\partial \dot{\bm q} } &= \dot{\bm u}^T + \bm u^T\frac{\partial \bm f}{\partial \dot{\bm q}} + \bm u^T \frac{\partial^2 \bm f}{\partial \dot{\bm q}^2}\dot{\bm q} - \bm u^T\frac{\partial \bm g}{\partial \dot{\bm q}} \label{dldqnew}
\end{align}

\noindent Taking the total derivative of Eq. \eqref{dldqnew} with respect to time yields 

\begin{equation}
    \begin{split}
        \frac{d}{dt}\left(\frac{\partial  L}{\partial \dot{\bm q}}\right) = \ddot{\bm u}^T &+ \dot{\bm u}^T\frac{\partial \bm f}{\partial \dot{\bm q}} +\bm  u^T\frac{d}{dt}\left(\frac{\partial\bm  f}{\partial \dot{\bm q}}\right) + \dot{\bm u}^T\frac{\partial^2 \bm f}{\partial \dot{\bm q}^2}\dot{\bm q}  \\ 
        &+ \bm u^T\frac{d}{dt}\left(\frac{\partial^2\bm  f}{\partial \dot{\bm q}^2}\right) \dot{\bm q}+ \bm u^T \frac{\partial^2 \bm f}{\partial \dot{\bm q}^2} \ddot{\bm q} \\
        &- \dot{\bm u}^T\frac{\partial \bm g}{\partial \dot{\bm q}} -\bm  u^T \frac{d}{dt}\left(\frac{\partial \bm g}{\partial \dot{\bm q}}\right)
    \end{split}
    \label{ddtdldqdot}
\end{equation}

\noindent By combining Eqs. \eqref{ddtdldqdot} and \eqref{dldqnew} and collecting terms, and utilizing \textbf{Theorem \ref{TheoremLagrange}} we arrive at the differential equations for the control from this variational principle. 

\begin{equation}
    \ddot{\bm u}^T = -\dot{\bm u}^T\frac{\partial^2 \bm f}{\partial \dot{\bm q}^2}\dot{\bm q} - \bm u^T\left[ \frac{d}{dt}\left(\frac{\partial^2 \bm f}{\partial \dot{\bm q}^2}\dot{\bm q}  \right) - \frac{\partial}{\partial \bm q}\frac{\partial \bm f}{\partial \dot{\bm q}}\dot{\bm q} + \frac{\partial \bm g}{\partial \bm q} - \frac{\partial \bm h}{\partial\bm q} \right]
    \label{uddoteqnlagrangian}
\end{equation}

\noindent As established above, $\bm g$ (and by extension $\bm f$) is quadratic in the generalized velocities, and so $\ds \frac{\partial \bm f}{\partial \dot{\bm q}}\dot{\bm q} = \frac{\partial \bm g}{\partial \dot{\bm q}}\partial \bm q =2\bm g$. Taking the derivative with respect to $\dot{\bm q}$, 
\begin{equation}
\begin{split}
    \frac{\partial}{\partial \dot{\bm q}} \left(\frac{\partial \bm f}{\partial \dot{\bm q} }\dot{\bm q}\right) &= \frac{\partial}{\partial \dot{\bm q}}(2\bm g) \\ 
    \frac{\partial^2 \bm f}{\partial \dot{\bm q}^2}\dot{\bm q} + \frac{\partial \bm f}{\partial \dot{\bm q}} &= 2\frac{\partial \bm g}{\partial \dot{\bm q}} \\ 
    \frac{\partial^2 \bm f}{\partial \dot{\bm q}^2}\dot{\bm q} &= \frac{\partial \bm f}{\partial \dot{\bm q}}
\end{split}
\label{d2fdq2}
\end{equation}

\noindent The relation $\ds \frac{\partial \bm f}{\partial \dot{\bm q}} = \frac{\partial \bm g}{\partial \dot{\bm q}}$ was used in the above equation. We can use Eq. \eqref{d2fdq2} to rewrite Eq. \eqref{uddoteqnlagrangian}. 

\begin{equation*}
\begin{split}\label{OCsolVar}
    \ddot{\bm u}^T &= -\dot{\bm u}^T \frac{\partial \bm f}{\partial \dot{\bm q}} - \bm u^T \left[ \frac{d}{dt}\left(\frac{\partial \bm f}{\partial \dot{\bm q}} \right) - \frac{\partial}{\partial \bm q}(2\bm g) + \frac{\partial \bm g}{\partial \bm q} - \frac{\partial \bm h}{\partial \bm q}\right] \\ 
    &= -\dot{\bm u}^T\frac{\partial \bm f}{\partial \dot{\bm q}} - \bm u^T\left[ \frac{d}{dt}\left(\frac{\partial \bm f}{\partial \dot{\bm q}} \right) - \frac{\partial \bm g}{\partial \bm q} - \frac{\partial \bm h}{\partial \bm q}\right] \\
    &= -\dot{\bm u}^T\frac{\partial \bm f}{\partial \dot{\bm q}} - \bm u^T\left[ \frac{d}{dt}\left(\frac{\partial\bm  f}{\partial \dot{\bm q}} \right) - \frac{\partial \bm f}{\partial \bm q} \right] \\
\end{split}
\end{equation*}

This form for the optimal control is identical to the one in Eq. \eqref{optimalControlODE}. The {stationary points} of $\int L\,dt$ {are therefore} the correct equations for the optimal control.

A similar argument is used to find the equations of motion for the states. Taking the derivative of Eq. \eqref{BetaFunctional} with respect to $\bm u $ and $\dot{\bm u}$ yields

\begin{align}
    \frac{\partial L}{\partial \bm u^T} &= \dot{\bm{q}} \\ 
    \frac{\partial L}{\partial \dot{\bm{u}}^T} &= \frac{\partial \bm f}{\partial \dot{\bm q}}\dot{\bm q} - \bm g +\bm  h + \bm{u} \label{dLduProof3}
\end{align}

\noindent We know from above that $\ds \frac{\partial \bm g}{\partial \dot{\bm q}}\dot{\bm q} = \ds 2\bm g$. Therefore, Eq. \eqref{dLduProof3} reduces to $2\bm g - \bm g + \bm h + \bm u = \bm g + \bm h + \bm u = \bm f + \bm u$. Combining everything together,

\begin{equation*}
    \frac{d}{dt}\left( \frac{\partial L }{\partial \dot{\bm u}} \right) - \frac{\partial L }{\partial \bm{u}} \implies \ddot{\bm{q}} = \bm f + \bm u
\end{equation*}

\noindent This is the correct equation of motion for the controlled system, defined in Eq. \eqref{SysDyn}. \setlength\fboxrule{1.1pt}\setlength\fboxsep{0pt}\fbox{\phantom{\rule{5pt}{5pt}}}    

The above proof used an extended version of the optimal control functional which explicitly shows the separation into the {position} dependent terms and velocity dependent terms to demonstrate how Noether's theorem is used to arrive at the functional. However, it is possible to simplify the form of the optimal control functional given in \eqref{BetaFunctional2} as described below. Define the operator
\begin{equation*}
    \psi_i(L) = \frac{d}{dt}\left(\frac{\partial L}{\partial \dot{q}_i} \right) - \frac{\partial L}{\partial q_i}
    \label{psieq}
\end{equation*}

\noindent Following Eq. \eqref{lagrangianEquation}, we can solve for the equations of motion of the uncontrolled system, $f$, in terms of the classical Lagrangian, $\mathcal{L} = T-V$.
\begin{align}
    \psi(\mathcal L) &= \alpha (\ddot{\bm q} - \bm f) \\ 
    \implies&\bm  f = \ddot{\bm q} - \alpha^{-1}\psi(\mathcal L)
    \label{fequationJLM}
\end{align}

Additionally, the condition in Eq. \eqref{quadraticSatisfaction} yields $$\frac{\partial \bm g}{\partial \dot{\bm q}}\dot{\bm q} = \frac{\partial \bm f}{\partial \dot{\bm q}}\dot{\bm q} = 2\bm g$$ Applying both this and the relation in Eq. \eqref{fequationJLM} to the functional $L$ in Eq. \eqref{BetaFunctional2} gives 

\begin{equation*}
    \begin{split}
        L &= \dot{\bm u}^T\dot{\bm q} + \bm u^T\frac{\partial \bm f}{\partial \dot{\bm q}} \dot{\bm q} - \bm u^T\bm g + \bm u^T \bm h + \frac{1}{2}\bm u^T\bm u \\ 
        &= \dot{\bm u}^T\dot{\bm q} + \bm u^T(2\bm g) - \bm u^T\bm g + \bm u^T\bm  h + \frac{1}{2}\bm u^T\bm u \\
        &= \dot{\bm u}^T\dot{\bm q} + \bm u^T\bm f + \frac{1}{2}\bm u^T\bm u \\
        &= \dot{\bm u}^T\dot{\bm q} + \bm u^T(\ddot{\bm q} - \alpha^{-1}\psi(L)) + \frac{1}{2}\bm u^T\bm u 
    \end{split}
\end{equation*}

We can now write the optimal control Lagrangian explicitly in terms of the classical Lagrangian, the control $\bm u$ and the generalized coordinates $\bm q$.

\begin{equation}
    L =  \dot{\bm u}^T\dot{\bm q} + \bm u^T(\ddot{\bm q} - \alpha^{-1}\psi(\mathcal L)) + \frac{1}{2}\bm u^T\bm u 
    \label{optimalControlLagrangianGeneral}
\end{equation}
The Lagrangian in Eq.~\eqref{optimalControlLagrangianGeneral} (or Eq.~\eqref{BetaFunctional2}) is in a general form that applies to any case where a control is applied to a conservative system, e.g. the thrust control on a spacecraft. Equation ~\eqref{optimalControlLagrangianGeneral} is used to get the optimal control by applying the same process as in Hamilton's Principle. We can do that in a general form as follows. First, we minimize the functional $\int L\,dt$ with respect to the generalized coordinates $\bm q$, by computing the partial derivatives of $L$ with respect to $\bm q$ and $\dot{\bm q}$. Substituting these into the Euler-Lagrange Equations will have us arrive at the differential equations for the optimal control in \eqref{optimalControlODE}.  

\subsection{Example}
To demonstrate the method, we consider the two body problem of gravitation. The problem describes the position of an orbiting body around a central mass in terms of polar coordinates $r,\theta$. This is a foundational problem in low thrust space trajectory design.



Suppose the desired control objective is to reach a final prescribed state while minimizing the control. The state vector is given by $\bm x^T = [r, \theta, \dot{r}, \dot{\theta}] = [\bm q^T, \dot{\bm q}^T].$ Formally, the optimal control problem is $$\mathrm{min}\,\, \mathcal{J} = \bm \nu^T\bm F + \int_{t_0}^{t_f}\frac{1}{2}(u_r^2 + u_{\theta}^2)\,dt$$ where $\bm{F}(t_f,\bm {x}(t_f))=0$ are the terminal state constraints.

To find the variational principle which will provide the optimal control, we must find the equations of motion and the Jacobi Last Multiplier matrix for the original uncontrolled system. The classical Lagrangian for this system and its Jacobi Last multiplier are given by

\begin{equation*}
    \mathcal L = \frac{1}{2}\left(\dot{r}^2 + r^2\dot{\theta}^2\right) + \frac{\mu}{r}
\end{equation*}
\begin{equation*}
    \alpha = \frac{\partial^2 \mathcal L}{\partial \dot{\bm q}^2} = \left[ \begin{array}{cc}
        1 & 0 \\
        0 & r^2
    \end{array}\right]
\end{equation*}
\noindent With this, we can substitute all values into Eq. \eqref{optimalControlLagrangianGeneral} to arrive at the new functional
\begin{equation*}
    \begin{split}
        L =  \dot u_r \dot r + \dot u_\theta \dot \theta+ \left(r\dot\theta^2 - \frac{\mu}{r^2}\right)u_r - \frac{2\dot r\dot \theta}{r}u_\theta+ \frac{1}{2}\left( u_r^2 + u_\theta^2\right)
    \end{split}
\end{equation*}

We then minimize this functional with respect to the generalized coordinates $q^T = [r,\theta]$ (following the variational approach briefed above) to arrive at the differential equations of optimal control for this system.

\begin{align}
    \ddot{u}_r &= \left(\dot\theta^2 + \frac{2\mu }{r^3}\right)u_r + \frac{2\ddot \theta}{r}u_\theta + \frac{2\dot\theta}{r}\dot u_{\theta} \label{urddot}\\ 
    \ddot{u}_{\theta} &= -(2\dot r \dot \theta + 2r\ddot \theta)u_r + \left( \frac{2\ddot r}{r} - \frac{2\dot r^2}{r^2}\right) u_{\theta}     - 2r\dot\theta\dot u_r + \frac{2\dot r}{r}\dot u_{\theta}
    \label{utddot}
\end{align}

{We can verify that this solution is the valid optimal control solution by comparing it against that returned by the classical optimal control theory. To begin, we start with the equations of motion of the spacecraft under a central gravitational field. 
\begin{align*}
    \ddot r &= r\dot\theta^2 - \frac{\mu}{r^2} +u_r\\ 
    \ddot \theta &= -\frac{2\dot r \dot \theta}{r} + u_\theta
\end{align*}
The optimal control Hamiltonian is given by $$H = \frac{1}{2}(u_r^2 + u_\theta^2) + \lambda_1\dot{r} + \lambda_2\dot{\theta} + \lambda_3\left(r\dot\theta^2 - \frac{\mu}{r}\right) - \frac{2\lambda_4\dot r\dot\theta}{r}$$ The costate equations of motion are found with Eq. \eqref{ELElambda}.
\begin{equation}
    \begin{split}
        \dot\lambda_1 &= -\lambda_3\left(\frac{2\mu}{r^3} + \dot\theta^2\right) - \frac{2\lambda_4\dot r\dot \theta}{r^2} \\ 
        \dot\lambda_2 &= 0 \\ 
        \dot\lambda_3 &= \frac{2\lambda_4\dot\theta}{r} - \lambda_1 \\ 
        \dot\lambda_4 &= \frac{2\lambda_4\dot r}{r} - 2\lambda_3 r \dot\theta - \lambda_2
    \end{split}
    \label{costatespacecraft}
\end{equation}
The optimal control law is found by applying the stationarity condition $\partial H/\partial\bm u=\bm0$, which gives $u_r=-\lambda_3$ and $u_\theta=-\lambda_4$.
It is desired to show that Eqs. \eqref{urddot}-\eqref{utddot} are equivalent to the costate differential equations in Eq. \eqref{costatespacecraft} with optimal control policy $u_r=-\lambda_3,u_\theta=-\lambda_4$. Taking the derivative of $\dot{\lambda}_3$ yields
\[
\ddot{\lambda}_3 = 2\dot{\lambda}_4\frac{\dot\theta}{r} + 2\lambda_4\left(\frac{\ddot\theta}{r} - \frac{\dot\theta}{r^2}\right)-\dot\lambda_1
\]
Substituting in the expression for $\dot\lambda_1$, and rewriting $\lambda_3$ and $\lambda_4$ end their derivatives in terms of $u_r,u_\theta,\dot{u}_r$, and $\dot u_\theta$ yields
\[
\ddot{u}_r = \left(\dot\theta^2 + \frac{2\mu }{r^3}\right)u_r + \frac{2\ddot \theta}{r}u_\theta + \frac{2\dot\theta}{r}\dot u_{\theta}
\]
This is the same equation of motion for $u_r$ as was derived using the generalized Lagrangian approach. We can take the derivative of $\dot\lambda_4$ to get 
\[
\ddot\lambda_4=\dot\lambda_4\left(\frac{2\dot r}{r}\right) + 2\lambda_4\left(\frac{\ddot r}{r} - \frac{\dot r^2}{r^2}\right) - 2\dot\lambda_3r\dot\theta - 2\lambda_3(\dot r\dot\theta + r\ddot\theta) - \dot\lambda_2
\]
Substituting in the equations for $\dot{\lambda}_2$, as well as the control variables and the derivatives yields
\[
\ddot{u}_{\theta} = -(2\dot r \dot \theta + 2r\ddot \theta)u_r + \left( \frac{2\ddot r}{r} - \frac{2\dot r^2}{r^2}\right) u_{\theta}     - 2r\dot\theta\dot u_r + \frac{2\dot r}{r}\dot u_{\theta}
\]
Once again, this is the same equation of motion for $u_\theta$ as was arrived at using the generalized Lagrangian. We note that the control equations of motion in Eqs. \eqref{urddot}-\eqref{utddot} did not require any application of optimal control theory; only Hamilton's law of varying action was required to reach these equations of motion. 
}

\section{Diagonalization of the Jacobi’s Last Multiplier Matrix (Inverse Lagrangian)}\label{SecInvLagrang}
Here we present another approach for deriving the Lagrangian of a system accounting for the control. For this, we will go into more detail into the Jacobi's Last Multiplier (JLM) matrix. Consider an arbitrary dynamical system which can be derived using a variational principle. As it is derivable from a variational principle, there exists an action functional $S =\int_{t_0}^{t_f} \mathcal B\,dt$ which, when minimized, returns the system equations of motion.  For an $n$ degree-of-freedom system, the JLM $\alpha$ is a symmetric $n\times n$ matrix defined as the Hessian of the functional's Lagrangian $\mathcal B$.
\begin{equation*}
    \alpha_{ij} = \frac{\partial^2 \mathcal B}{\partial \dot q_i \partial \dot q_j},\,\,\,\,i,j=1,2,...,n
    \label{JLM}
\end{equation*}
where $q_i$, for $ i= 1, \cdots, n,$ are the generalized coordinates.

In the field of analytical mechanics, the inverse Lagrangian problem is that given a set of $n$ differential equations $\ddot{q}_i = f_i(t,\bm{q},\dot{\bm{q}})$, find a Lagrangian which, when rendered stationary, yields this set of differential equations.
In fact, the literature is generally mainly concerned with proving the existence of the Lagrangian, as finding the Lagrangian itself can be done through an over determined set of partial differential equations resulting from the Euler-Lagrange equations of the system. Here, the focus is different; the goal is to write a modified Lagrangian for the system that accounts for states and control. 
The Lagrangian in the above variational principle will not be the traditional Lagrangian. In this sense, this proposed approach can be seen as the search for a {nonstandard Lagrangian} which describes an optimal control system. 

One question which arises when searching for a nonstandard Lagrangian is that of knowing if the Lagrangian is correct. The criterion for a Lagrangian to be correct is that the correct equations of motion arise when the Lagrangian is passed through Hamilton's law of varying action (or any variational principle) \cite{sarlet}. One can also verify that a Lagrangian is correct by verifying the Helmholtz conditions \cite{HelmHotz1941Douglas} on the dynamic system, which are here described. Suppose that an $n$ degree of freedom system is governed by the set of second order differential equations $\ddot q_i = f_i(t,\bm{q},\dot{\bm{q}})$. The Helmholtz conditions are a set of necessary and sufficient conditions such that a nonsingular multiplier matrix $\alpha$ exists such that 
\begin{equation}
    \alpha_{ij}(\ddot{q}_j-f_j) = \frac{d}{dt}\left(\frac{\partial L}{\partial \dot{q}_i}\right) - \frac{\partial L}{\partial q_i} 
    \label{HelmholtLagrange}
\end{equation}
 As described above, the matrix $\alpha$ is the Hessian of the Lagrangian with respect to the generalized velocities. The Helmholtz conditions are given in \cite{sarlet} as
\begin{equation}
    \alpha_{ij} = \alpha_{ji} 
    \label{helmholtzConditions1}
\end{equation}
\begin{equation}
        \frac{\partial \alpha_{ij}}{\partial \dot{q}_k} = \frac{\partial \alpha_{ik}}{\partial \dot{q}_j} 
    \label{helmholtzConditions2}
\end{equation}
\begin{equation}
        \frac{\partial \beta_i}{\partial \dot{q}_j} + \frac{\partial \beta_j}{\partial \dot{q}_i} = 2\frac{d}{dt}(\alpha_{ij}) 
    \label{helmholtzConditions3}
\end{equation}
\begin{equation}
\frac{\partial \beta_i}{\partial q_j} + \frac{\partial \beta_j}{\partial q_i} = \frac{1}{2}\frac{d}{dt}\left(\frac{\partial \beta_i}{\partial \dot{q}_j} - \frac{\partial \beta_j}{\partial \dot{q}_i}\right)
    \label{helmholtzConditions4}
\end{equation}

\noindent where $\beta_i = -\sum_{j}\alpha_{ij}f_j$. These conditions must be satisfied for all combinations of $i,j,k\in \{1,2,...,n\}$. If the correct equations of motion are returned, and the Jacobi's Last Multiplier matrix (along with the equations of motion) satisfy the Helmholtz conditions, then the Lagrangian is a correct representation of the system's dynamics. Note that Lagrangians will not be unique, and a single problem which admits one Lagrangian will admit infinitely many. If $L$ describes the dynamics of a system, then so will $L'=L + \dot g$, where $g(t,\bm{q},\bm{\dot{q}})$ is an arbitrary differentiable function \cite{neuenschwander}. This detail will become relevant when searching for new invariances on the optimal path.

For typical physical dynamic systems, the classical Lagrangian is the excess of kinetic energy over potential energy. This, in turn, gives a leading velocity term that is quadratic in the generalized velocities, $L = \frac{1}{2}\sum\alpha_{ii}\dot{q}_i^2 + ...$. In such a Lagrangian, each of the generalized velocities is only multiplied by itself in the Lagrangian. This leads to a JLM matrix which has nonzero entries in the main diagonal, and zeroes elsewhere. {In this proposed approach}, the JLM would instead have zeroes on the main diagonal, and nonzero entries elsewhere in the matrix. This is indicative of the direct coupling of the state and control variables in the Lagrangian. If in the new system ($F'$ in Fig. ~\ref{BlockDiag}) $\hat{\bm{q}}=[\bm{q}^T, \bm{u}^T]^T$, then the Lagrangian would take the form $L=\sum_{i\neq j}\alpha_{ij}\dot{\hat{q}}_i \dot{\hat{q}}_j' + ...$. The only nonzero elements of the JLM are those that contain a state and a control generalized velocity. On the upper triangular part of the matrix, this means that $\alpha_{ij}=0$ always (but not exclusively) when $i \leq j$.

When computing the equations of motion for this system, treating both the original states and the controls as a new state vector $\hat{\bm{q}}$, the equations of motion for $\bm{q}$ arise when taking the variation with respect to $\bm{u}$, and the equations for the control variables arise when taking the variation with respect to the original state variables, $\bm{q}$. We can arrive at conditions on the Lagrangian by applying Hamilton's law of varying action to this new action functional. $$\delta S' = \int_{t_0}^{t_f}\left( \frac{\partial L}{\partial \bm{q}}\delta\bm{q} + \frac{\partial L }{\partial \dot{\bm{q}}}\delta\dot{\bm{q}} + \frac{\partial L}{\partial \bm{u}}\delta\bm{u} + \frac{\partial L }{\partial \dot{\bm{u}}}\delta\dot{\bm{u}} \right)\,dt =0$$ Integrating by parts on the $\delta\dot{\bm{q}}$ and $\delta\dot{\bm{u}}$ terms, 
\begin{align*}
    \delta S' = \int_{t_0}^{t_f}\left( \,\, \left\{\frac{d}{dt}\left(\frac{\partial L}{\partial \dot{\bm{q}}}\right)  - \frac{\partial L}{\partial \bm{q}} \right\}\delta\bm{q} \right.+   
    \left.\left\{ \frac{d}{dt}\left(\frac{\partial L}{\partial \dot{\bm{u}}}\right) - \frac{\partial L}{\partial \bm{u}} \right\}\delta\bm{u}\right)& = 0
\end{align*}

\noindent The variations in $\bm{u}$ and $\bm{q}$ are arbitrary and unrelated, so we get the well known Lagrangian formulation for this problem. From the structure of the JLM that was asserted above, we know that the variation with respect to the control will yield the state equations of motion, and vice versa. From equation \eqref{HelmholtLagrange}, we know what the state equations of motion have to be. Therefore, we arrive at the condition that the Lagrangian must satisfy

\begin{align*}
    \frac{d}{dt}\left(\frac{\partial L}{\partial \dot{u_i}}\right) - \frac{\partial L}{\partial u_i} &=  \frac{\partial V}{\partial q_j} + \frac{d}{dt}\left(\frac{\partial T}{\partial\dot{q_j}}   \right) -  \frac{\partial W_{nc}}{\partial q_j} 
\end{align*}

The values of $i$ and $j$ which make the equations of motion valid are dependent on the JLM matrix. Specifically, they will correspond to the entries $\alpha_{ij}$ which are nonzero.

We will not know in general ahead of time where these nonzero entries in the JLM will be, nor what their value will be. Instead, when using this inverse Lagrangian approach, we search for an invertible coordinate transformation on the equations of motion $f_i(t,\bm{q},\bm{\dot{q}})\mapsto\Tilde{f}_i(t,\Tilde{\bm{q}}, \dot{\Tilde{\bm{q}}})$ such that the identity matrix becomes a Jacobi Last Multiplier for the transformed system. With the identity matrix satisfying the Helmholtz conditions, the Lagrangian can then be found in the transformed coordinates.
The Lagrangian in the transformed coordinates can then be expressed as 
\begin{equation}
    \begin{split}
\Tilde{{L}} = \frac12\sum\dot{\Tilde{q}}_i^2 + \sum\dot{\Tilde{q}}_i g(t,\Tilde{\bm{q}}) + h(t,\Tilde{\bm{q}})
    \end{split}
    \label{ModLagr}
\end{equation}
Where the $g(t,\Tilde{\bm{q}})$ and $h(t,\Tilde{\bm{q}})$ functions show up when carrying out the integration. This Lagrangian can now be transformed back to the original coordinates through the inverse of the original coordinate transformation. This inverse coordinate transformation, however, is problem dependent because it is dependent on the equations of motion. We will demonstrate the coordinate transformation approach via an example.

\subsection*{Example} 
\begin{figure}[hbt!]
    \centering
    \includegraphics[width=0.5\linewidth]{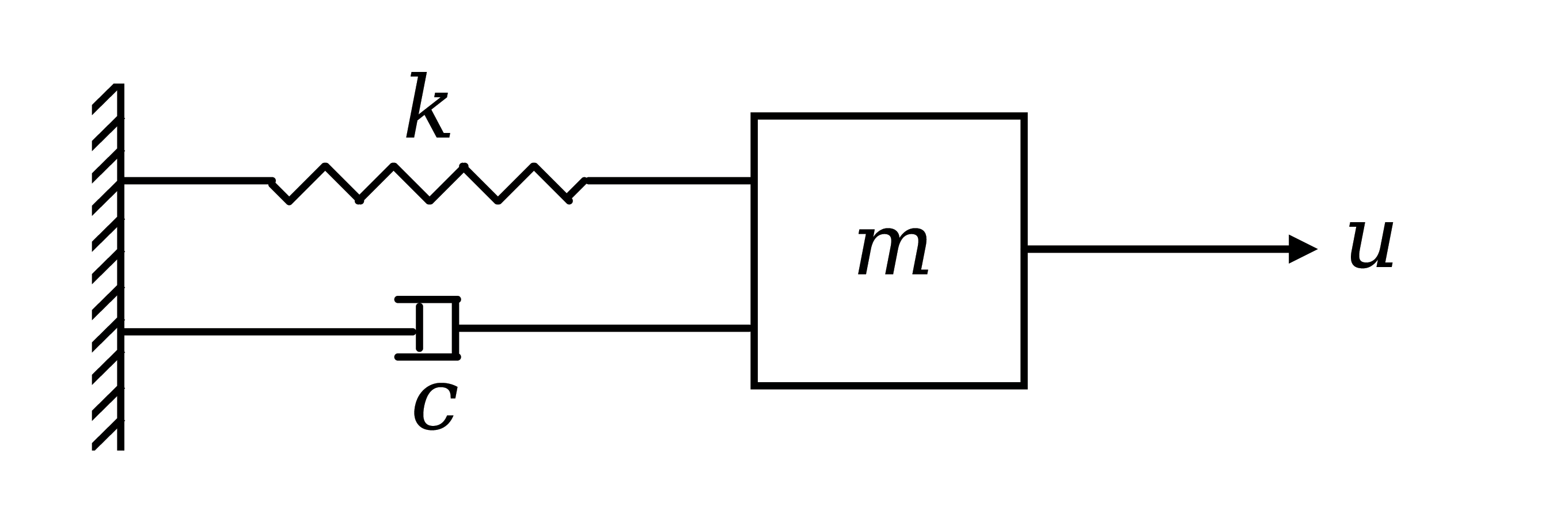}
    \caption{Mass-spring-damper system}
    \label{fig:mckDiagram}
\end{figure}

Consider a mass-spring-damper system shown in Figure \ref{fig:mckDiagram}. The system has a spring constant $k$, damping coefficient $c$, mass $m$, and a control force $u$ applied to the mass. The position of the mass is represented by $x$. 
Suppose the desired control objective is to reach a final prescribed state while minimizing the control. Formally, $$\mathrm{min}\,\, \mathcal{J} = \bm\nu^T\bm F + \int_{t_0}^{t_f}\frac{1}{2}u^2\,dt$$ where $\bm F(t_f,\bm{x}(t_f))=0$ are the terminal state constraints. The governing system equations of motion are $m\ddot{x} = c\dot{x} +kx = u$. Here we assume we know the equation of optimal control, and use it to find the Lagrangian. It is possible to use classical optimal control theory to show that for this example the optimal control must satisfy: 
\begin{equation*}
    m\ddot{u} - c\dot u + ku = 0
\end{equation*}

We now search for a coordinate transformation $\bm y = T\bm q$ which will give us an identity Jacobi Last Multiplier matrix in the transformed coordinates $\bm  y$. Let 

\begin{equation}
    T = \left[ \begin{array}{cc} f_1(t,x,u)\,\,\,\,\, & g_1(t,x,u) \\ f_2(t,x,u)\,\,\,\,\, & g_2(t,x,u) \end{array} \right],\,\,\,\,\,\, T^{-1} = \frac{1}{f_1g_2-f_2g_1}\left[ \begin{array}{cc} g_2 & -g_1 \\ -f_2 & f_1 \end{array} \right] = \left[ \begin{array}{cc} K & L \\ M & N \end{array}\right]
    \label{coordTransform}
\end{equation}

\noindent where $\bm q = [x,u]^T$. For this to be an invertible transformation we require that $D \equiv f_1g_2-f_2g_1\neq 0$. By the product rule, 

\begin{equation*}
    \begin{split}
        \ddot{y}_1 &= \ddot{f}_1x + 2\dot{f}_1\dot{x} + f_1\ddot{x} + \ddot{g}_1u+2\dot{g}_1\dot{u} + g_1\ddot{u} \\ 
        \ddot{y}_2 &= \ddot{f}_2x + 2\dot{f}_2\dot{x} + f_2\ddot{x} + \ddot{g}_2u+2\dot{g}_2\dot{u} + g_2\ddot{u}
    \end{split}
\end{equation*}

\noindent Substituting in the differential equations, as well as the inverse coordinate transformation in Eq. \eqref{coordTransform}, we arrive at the equations of motion in the new coordinates.

\begin{equation*}
    \begin{split}
        \ddot{y}_1 &= [\ddot{f}_1-\omega^2f_1](Ky_1+Ly_2) + [2\dot{f}_1-z^2f_1](\dot{K}y_1+K\dot{y}_1 + \dot{L}y_2 + L\dot{y}_2) \\ 
        &+\left[\frac{f_1}{m}+\ddot{g}_1-\omega^2 g_1\right](My_1 + Ny_2) + [2\dot{g}_1 + z^2g_1](\dot{M}y_1 + M\dot{y}_1 + \dot{N} y_2 + N\dot{y}_2) \\
        \ddot{y}_2 &= [\ddot{f}_2-\omega^2f_2](Ky_1+Ly_2) + [2\dot{f}_2-z^2f_2](\dot{K}y_1+K\dot{y}_1 + \dot{L}y_2 + L\dot{y}_2) \\ 
        &+\left[\frac{f_2}{m}+\ddot{g}_2-\omega^2 g_2\right](My_1 + Ny_2) + [2\dot{g}_2 + z^2g_2](\dot{M}y_1 + M\dot{y}_1 + \dot{N} y_2 + N\dot{y}_2)
    \end{split}
\end{equation*}

\noindent where $\omega^2=k/m, z^2=c/m$. In the new $(y_1,y_2)$ coordinates, we are now looking for values of $f_1, f_2, g_1, g_2$ which will make $\alpha = I$ be a valid multiplier. To do this, we apply the Helmholtz conditions in Eqs. \eqref{helmholtzConditions1} - \eqref{helmholtzConditions4}. Due to the choice of multiplier matrix, some of the conditions are trivially satisfied. The remaining conditions are 

\begin{equation*}
    \begin{split}
        \frac{\partial \beta_1}{\partial \dot{y}_1} &= 0 \\ 
        \frac{\partial \beta_2}{\partial \dot{y}_2} &= 0 \\ 
        \frac{\partial \beta_1}{\partial \dot{y}_2} + \frac{\partial \beta_2}{\partial \dot{y}_2} &= 0 \\
        \frac{\partial \beta_1}{\partial y_2} - \frac{\partial \beta_2}{\partial y_1} &= \frac{1}{2}\frac{d}{dt}\left[ \frac{\partial \beta_1}{\partial \dot{y}_2} - \frac{\partial \beta_2}{\partial \dot{y}_1}\right]
    \end{split}
\end{equation*}

\noindent Adding together the first two conditions yields $\dot{f}_1g_2 + f_1\dot{g}_2 - \dot{f}_2g_1 - f_2\dot{g}_1 =\frac{d}{dt}(f_1g_2-f_2g_1)= 0$. We can recognize this as the total time derivative of $D$. The first two conditions together imply that $D\equiv\mathrm{const.}$ We also know from the requirement that $T$ be an invertible mapping that $D$ must be nonzero. With the knowledge that $D$ is a nonzero constant, after carrying out the necessary algebra the last condition can compactly be expressed as $f_1^2 + f_2^2 = 0$. Finally, the third condition requires that $f_1\dot{g}_1-\dot{f}_1g_1 + z^2f_1g_1 + \dot{f}_2g_2-f_2\dot{g}_2 - z^2f_2g_2=0$. These conditions will all be satisfied by a coordinate transformation 

\begin{equation*}
    T = \left[ \begin{array}{cc} e^{-z^2t} & ie^{z^2t} \\ -ie^{-z^2t} & -e^{-z^2t} \end{array}\right]
\end{equation*}

\noindent where $i\equiv \sqrt{-1}$ is the base imaginary unit. Explicitly written out, this gives the coordinates 
\begin{equation}
    \begin{split}
        y_1 = e^{-\frac{ct}{m}}x + ie^{\frac{ct}{m}}u \, \, , \, \,
        y_2 = -ie^{-\frac{ct}{m}}x - e^{\frac{ct}{m}}u
    \end{split}
    \label{transfomration}
\end{equation}

\noindent Under these transformed coordinates, the equations of motion are given by 
\begin{equation}
    \begin{split}
        \ddot y_1 &= z^2 i\dot{y}_2 + \left(\omega^2 - \frac{i}{2m}\right)y_1 - \frac{1}{2m}y_2 \, \, , \, \, \\
        \ddot y_2 &= -z^2i\dot{y}_1 - \frac{1}{2m}y_1 + \left( \omega^2 + \frac{i}{2m}\right) y_2
    \end{split}
    \label{yddot}
\end{equation}
\noindent As $\alpha = I$ in the transformed coordinates, the Lagrangian in the new coordinates is of the form given in Eq.~\eqref{ModLagr}: $$\Tilde{{L}} = \frac12\dot{y}_1^2 + \frac12\dot{y}_2^2 + h_1(t,y_1,y_2,\dot{y}_2)\dot{y}_1 + h_2(t,y_1,y_2,\dot{y}_1)\dot{y}_2 + p(t,y_1,y_2)$$ Applying Eq. \eqref{HelmholtLagrange} to this Lagrangian and requiring that it return the equations of motion in Eq. \eqref{yddot} yields 
\begin{equation*}
    \tilde{{L}} = \frac{1}{2}\dot{y}_1^2 + \frac{1}{2}\dot{y}_2^2 - \frac{1}{2}z^2iy_2\dot{y}_1 + \frac{1}{2}z^2iy_1\dot{y}_2 - \frac{1}{2}\left(\frac{i}{2m}-\omega^2\right)y_1^2 - \frac{y_1y_2}{2m} + \frac{1}{2}\left(\omega^2 + \frac{i}{2m}\right)y_2^2
\end{equation*}

\noindent Applying the inverse transformation in Eq. \eqref{transfomration}, substituting back the values for $z$ and $\omega$, and dividing through by $2i/m$ yields the Lagrangian in the original coordinates:
\begin{equation}
    L = m\dot x\dot u - \frac12c\dot x u + \frac12cx\dot u - kxu + \frac12u^2
    \label{mcklagrangian}
\end{equation}
Hence the functional to be minimized is $$S'=\int_{t_0}^{t_f}\left(m\dot x \dot u - \frac{1}{2}c\dot x u + \frac{1}{2}c x \dot u - kxu + \frac{1}{2}u^2\right)\,dt$$ As described above, both the position $x$ and the control input $u$ are now considered generalized coordinates for the system. The integral $S'$ is the modified action functional which describes the dynamics of the system. Following the process of Hamilton's law of varying action, we take the complete variation of the action functional and assert that it vanishes to minimize the functional. 
\[
\begin{split}
    \delta S' = \int_{t_0}^{t_f}\left\{ \left(m\dot u - \frac{1}{2}cu\right) \right. \delta\dot x + \left(m\dot x + \frac{1}{2}cx\right) \delta \dot u  \,\,\,\,& \\
    + \left.\left(\frac{1}{2}c\dot u - ku \right)\delta x + \left(u-\frac{1}{2}c\dot x - kx \right)\delta u  \right\}& \,dt=0
\end{split}
\]

After carrying out the necessary integration by parts on the $\delta \dot x$ and $\delta \dot u$ terms, the total variation of the functional can be written as 

\[
\begin{split}
    \delta S' &= \left\{ \left( m\dot u - \frac{1}{2}cu\right)\delta u + \left( m\dot x + \frac{1}{2}cx\right)\delta u \right\}_{t_0}^{t_f} \\
    &+\int_{t_0}^{t_f} \left \{\left(-m\ddot u + c\dot u - ku\right) \delta x + \left(-m\ddot x - c\dot x - kx + u \right) \delta u \right \} \, dt =0
\end{split}
\]

Under the assumptions of Hamilton's law of varying actions, the variations are chosen such that $\delta q_i(t_0)=\delta q_i(t_f) = 0$ for all $1\leq i \leq n$. As such, the boundary terms in the variation of the action functional vanish. This leaves us with the two bracketed terms. In order for the integral to vanish for any arbitrary variation, the bracketed terms must also identically vanish. This gives us the equations of motion for our expanded generalized coordinates $\bm{\hat{q}}=[x,u]^T$.
\begin{equation}\label{SMDEqnMotion}
    m\ddot x + c\dot x + kx = u
\end{equation}
\begin{equation}\label{SMDOptCtrl}
    m\ddot u - c\dot u + ku = 0
\end{equation}
Eq.~\eqref{SMDEqnMotion} is the typical equation of motion for a spring-mass-damper system. Eq.~\eqref{SMDOptCtrl}  describes the time evolution of the control which satisfies the necessary conditions for optimality and minimizes the optimal control objective $\mathcal{J}$. 
{With this approach we are able to write an equation for the optimal control without having to use the dynamic costates. } 

\section{Invariances on the Optimal Path}\label{invariancesSection}

Once a Lagrangian is obtained for an optimal control problem, it can be used to arrive at new conserved quantities for the optimal control system in terms of the state and control variables of the optimal control system. To find these invariances, we introduce the generator of the infinitesimal transformation in Eqs. \eqref{NoetherThermTimeTrans} and \eqref{NoetherThermCoordTrans}. The $(t,\bm{q})$ variables will be transformed into the $(t',\bm{q}')$ variables by the \textit{infinitesimal} operator

\begin{equation}
    G(t,\bm{q}) = \xi(t,\bm{q})\frac{\partial}{\partial t} + \sum_{i=1}^n \eta_i(t,\bm{q}) \frac{\partial }{\partial q_i}
    \label{infinitesimalGenerator}
\end{equation}

The original time and coordinate variables $(t,\bm{q})$ are transformed to the new variables $(t',\bm{q}')$ through the one parameter group of \textit{finite} transformations

\begin{align*}
    t' = e^{\varepsilon G(t,\bm{q})}t \\ 
    q_i' = e^{\varepsilon G(t,\bm{q})}q_i
\end{align*}

\noindent Note that the \textit{infinitesimal} transformation given in Eqs. \eqref{NoetherThermTimeTrans} and \eqref{NoetherThermCoordTrans} is the linear term of the Taylor expansion of the \textit{finite} transformations above around $\varepsilon=0$. To find invariances, we must find the transformations which leave the action functional unchanged under the action of the group. For this, we introduce the generator of the first extended group.

\begin{equation*}
    E(t,\bm{q},\dot{\bm{q}}) = \xi(t,\bm{q})\frac{\partial}{\partial t} + \sum_{i=1}^n \eta_i(t,\bm{q}) \frac{\partial }{\partial q_i} + \left[ \dot \eta_i - \dot\xi\dot q_i\right]\frac{\partial }{\partial \dot{q}_i}
    \label{extendedgroupgenerator}
\end{equation*}

\noindent where $\ds \dot\xi = \ds \frac{\partial \xi}{\partial t}+\sum_i \frac{\partial \xi }{\partial q_i}\dot{q}_i$ and $\ds \dot\eta_i =\ds\frac{\partial \eta_i}{\partial t} + \sum_m \frac{\partial \eta_i}{\partial q_m}\dot{q}_m$. The development and interpretation of these generators can be found in Ref. \cite{cohen1911introduction}. Consider the function 

\begin{equation}
    \Phi = \frac{\partial L}{\partial \dot{\bm{q}}}(\xi \dot{\bm{q}} - \bm{\eta}) - \xi L + \varphi(t,\bm{q})
    \label{noetherinvariance}
\end{equation}

\noindent Noether's theorem states that $\Phi$ will be a conserved quantity when the action integral defined by $L$ is both stationary and invariant. This form with the additional function $\varphi(t,\bm{q})$ is known as the ``divergence invariant" form \cite{neuenschwander} of Noether's theorem. It can be shown \cite{lutzky1978symmetry} that $\Phi$ will be an invariant when $\xi,\bm{\eta},$ and $\varphi$ can be chosen such that 

\begin{equation}
    E\{L\} = -\dot{\xi}L + \dot \varphi
    \label{invarianceIdentity}
\end{equation}

Equation \eqref{invarianceIdentity} is a scalar equation, yet it will provide a system of partial differential equations (PDEs) in $\xi,\bm{\eta},$ and $\varphi$. This scalar equation provides a system of PDEs because each velocity dependent term in the equation must separately vanish. This is because, following from Noether's theorem, we are only allowing our infinitesimal transformations to be position dependent, not velocity dependent. Once these PDEs are solved, each solution will have an associated conservation law from Noether's theorem. We proceed with an example.

Suppose we have a controlled mass-spring system (Fig. \ref{fig:mckDiagram} with a damping coefficient of $c=0$.) Then, from either Eq. \eqref{mcklagrangian} or  Eq. \eqref{optimalControlLagrangianGeneral}, we can write the Lagrangian for the optimal control system as $$L = m\dot{x}\dot{u} - kxu + \frac12u^2$$ The extended group generator for this system is 

\begin{equation*}
    E=\xi\frac{\partial }{\partial t} + \eta_1\frac{\partial }{\partial x} + \eta_2 \frac{\partial}{\partial u} + [\dot{\eta}_1 - \dot \xi \dot x]\frac{\partial }{\partial \dot x} + [\dot{\eta}_2 - \dot{\xi}\dot{u}]\frac{\partial}{\partial \dot{u}}
\end{equation*}

\noindent where we are using $\bm{q} = [x, u]^T$ as the complete set of generalized coordinates. Applying the condition in Eq. \eqref{invarianceIdentity} to the Lagrangian yields

\begin{equation*}
\begin{split}
    E\{L\} &= \eta_1(-ku) + \eta_2(u-kx) + [\dot{\eta}_1 - \dot \xi \dot x](m\dot u) + [\dot{\eta}_2 - \dot{\xi}\dot{u}](m\dot x) \\
    &= \dot{\xi}(m\dot x \dot u - k x u + \frac{1}{2}u^2) + \dot{\varphi}
\end{split}
\end{equation*}

We can expand the total time derivatives in this equation in terms of partial derivatives and collect terms to arrive at the following polynomial in generalized velocities. The notation $a_b\equiv\partial a/\partial b$ is used to represent the partial derivatives in the equation.

\begin{equation*}
    \begin{split}
        &(-ku\eta_1 + (u-kx)\eta_2) + \dot u[m\eta_{1,t}] + \dot{u}^2[m\eta_{1,u}] + \dot{x}[m\eta_{2,t}] + \dot{x}^2[m\eta_{2,x}]   \\ 
        &\,\,\,\,\,\,\,\,+\dot{x}\dot{u}[m\eta_{1,x} - 2m\xi_t + m\eta_{2,u}] + \dot{x}^2\dot{u}[-2m\xi_x] + \dot{x}\dot{u}^2[-2m\xi_u] \\ 
        &= (kxu\xi_t-\frac{1}{2}u^2\xi_t + \varphi_t) + \dot{u}[kxu \xi_u - \frac12u^2\xi_u + \varphi_u] + \dot{x}[kxu\xi_x - \frac12u^2\xi_x + \varphi_x]  \\
        &\,\,\,\,\,\,\,\,+\dot x \dot u [-m\xi_t] + \dot{x}^2\dot{u}[-m\xi_x] + \dot{x}\dot{u}^2[-m\xi_u^2]
    \end{split}
    \label{killingEquations}
\end{equation*}

\noindent As the condition must hold for all generalized velocities, the coefficients of each separate generalized velocity product term (including the constant term) must separately vanish. This leads to the following set of differential equations, known as the Killing equations. 

\begin{equation*}
    \begin{split}
        -ku\eta_1 + (u-kx)\eta_2 &= kxu\xi_t-\frac{1}{2}u^2\xi_t + \varphi_t \\ 
        m\eta_{2,t} &= kxu\xi_x - \frac12u^2\xi_x + \varphi_x \\
        m\eta_{1,t} &= kxu \xi_u - \frac12u^2\xi_u + \varphi_u \\
        m\eta_{2,x} &= 0 \\
        m\eta_{1,u} &= 0 \\ 
        m\eta_{1,x} - 2m\xi_t + m\eta_{2,u} &= -m\xi_t \\
        -2m\xi_x &= -m\xi_x \\
        -2m\xi_u &= -m\xi_u^2 \\
    \end{split}
\end{equation*}

\noindent It can be shown that this set of partial differential equations has the general solution 

\begin{equation*}
    \begin{split}
        \eta_1 &= c_1\left( \frac{\cos{\omega t} + 2\omega t\sin{\omega t}}{4\omega^2} \right) + c_2 \left( \frac{\sin{\omega t} - 2\omega t \cos{\omega t}}{4\omega^2}\right) \\ 
        &+ c_3\cos{\omega t} - c_4\sin{\omega t}  \\ 
        \eta_2 &= c_1\cos{\omega t} + c_2 \sin{\omega t }\\
        \xi&= c_5 \\ 
        \varphi &= c_1\left[\left( \frac{u-kx}{\omega} - \frac{3ku}{4\omega^3}\right)\sin{\omega t} + \frac{ku}{2\omega^2}\cdot t\cos{\omega t}\right] \\ 
        &+ c_2\left[ \left( \frac{kx-u}{\omega} + \frac{3ku}{4\omega^3}\right)\cos{\omega t} + \frac{ku}{2\omega^2}\cdot t\sin{\omega t}\right] \\ 
        &+ c_3\left[ \frac{-ku}{\omega}\sin{\omega t}\right] + c_4 \left[ \frac{-ku}{\omega}\cos{\omega t}\right]
    \end{split}
\end{equation*}

\noindent where $\omega \equiv \sqrt{k/m}$, and $c_1-c_5$ are arbitrary constants of integration. We can construct five linearly independent group generators from this solution, by letting in turn each $c_i=1$ and all others vanish. Each group generator will have an associated conservation law. These group generators are all of the form \eqref{infinitesimalGenerator} and are given by 

\begin{align*}
    G_1 &= \left( \frac{\cos{\omega t} + 2\omega t \sin{\omega t}}{4\omega^2}\right)\frac{\partial }{\partial x} + (\cos{\omega t})\frac{\partial }{\partial u} \\ 
    G_2 &= \left( \frac{\sin{\omega t} - 2\omega t\cos{\omega t}}{4\omega^2}\right) \frac{\partial }{\partial x} + (\sin{\omega t})\frac{\partial }{\partial u} \\ 
    G_3 &= (\cos{\omega t})\frac{\partial}{\partial x} \\ 
    G_4 &= (-\sin{\omega t})\frac{\partial}{\partial x} \\ 
    G_5 &= \frac{\partial}{\partial t}
\end{align*}

\noindent To each one of these group generators is associated a corresponding conserved quantity. By applying Eq. \eqref{noetherinvariance} to each generator, we can arrive at these invariances straightforwardly. For example, the conserved quantity from $G_5$ is simply the total energy-like term $\Phi_5=m\dot{x}\dot{u} + kxu - \frac{1}{2}u^2$ which we used in section \ref{InvarHamilt} to arrive at the optimal control Lagrangian in the energy case. For $G_1$, the associated conserved quantity is 

\begin{equation*}
    \Phi_1 = \left(\frac{u-kx}{\omega} - \frac{3ku}{\omega^3} - \frac{m\dot{u}}{2\omega}t\right)\sin{\omega t} + \left(\frac{ku}{2\omega^2}t - \frac{m\dot{u}}{4\omega^2} - m\dot{x}\right)\cos{\omega t}
\end{equation*}

\noindent Similarly, from $G_3$, we have a conserved quantity

\begin{equation*}
    \Phi_3 = -u\sqrt{mk}\sin{\omega t} - m\dot{u} \cos{\omega t}
\end{equation*}

\noindent The set $\Phi_1-\Phi_5$ is a set of functions of the state and control trajectories which are constant only for the optimally controlled system.

\section{Conclusions}
This paper presents a new approach developed to solve optimal control problems using variational methods of analytical mechanics. The new process starts by considering the control as a state of the system, and accounting for the control energy when writing the Lagrangian for the new system. 
Two methods are developed to derive the Lagrangian of the system accounting for the control. The first method leverages the fact that the Hamiltonian is invariant on the optimal trajectory and finds a Lagrangian function whose system equations of motion render the Hamiltonian invariant, within the context of Noether’s theorem. 
It is highlighted that this approach eliminates the need to compute dynamic costates in finding the optimal control. 
This paper also showed a method to write the Lagrangian of the system when accounting for the control; when the optimal control is known. This method finds a transformation in which the Jacobi’s Last Multiplier is an identity matrix, and use that to derive the Lagrangian in the original coordinate system. Finally, the paper shows how it is possible to use this Lagrangian to derive conserved quantities on the optimal path with Noether's Theorem. 

\section*{Data Availability Statement}
Data generated during the current study are available from the corresponding author on reasonable request. 

\section*{Acknowledgments}
This material is based upon work supported by the National Science Foundation Graduate Research Fellowship Program under Grant No. 2336877.

\bibliographystyle{unsrt} 
\bibliography{ORefs.bib,OCP.bib,Optimization,myRefs}

\end{document}